\documentclass[a4paper,10pt]{article}
\usepackage{amsfonts, amsthm}
\newcommand{\field}[1]{\mathbb{#1}}

\title{On analytic descriptions of two-dimensional surfaces associated 
with the $\field{C}P^{N-1}$ sigma model}

\author{\vspace{1cm}\\
         {\bf A.\ M.\ Grundland}$^{1,2}$
          \thanks{E-mail address:
           grundlan@crm.umontreal.ca} 
         {and \bf \.{I}.\ Yurdu\c{s}en}$^1$
             \thanks{E-mail address:
       yurdusen@crm.umontreal.ca}
          \\
          \\$^1$Centre de Recherches Math\'{e}matiques, Universit\'{e} 
                   de Montr\'{e}al, 
               \\ CP 6128, Succ.\ Centre-Ville, Montr\'{e}al, 
                        Qu\'{e}bec H3C 3J7, Canada
                   \\
                 \\$^2$Universit\'{e} du Qu\'{e}bec, 
                     Trois-Rivi\`{e}res, CP500, QC, G9A 5H7, Canada}

\date{\today}
\begin{document}

\maketitle

\begin{abstract}
We study analytic descriptions of conformal immersions of the 
Riemann sphere $S^2$ into the $\field{C}P^{N-1}$ sigma model. In 
particular, an explicit expression for two-dimensional ($2$-D) surfaces, 
obtained from the generalized Weierstrass formula, is given. It is also 
demonstrated that these surfaces coincide with the ones obtained from 
the Sym-Tafel formula. These two approaches correspond to parametrizations 
of one and the same surface in $\field{R}^{N^2-1}$.
\end{abstract}



%\newpage

\vskip 0.9cm
%\section{Formulation of the problem \label{formulation}}
In this note we investigate the relations between the 
$\field{C}P^{N-1}$ sigma model and the generalized Weierstrass 
formula for immersion of $2$-D surfaces in multi-dimensional 
Euclidean spaces. These links have been discussed in 
\cite{ghy, gy, Grundland2} and are governed by the formula 
\begin{eqnarray}
X_k(\xi, \bar{\xi})=
i\int_{\gamma} (-[\partial P_k, P_k]d\xi + [\bar{\partial}P_k, P_k] d\bar{\xi})\,,
\qquad k=0,1,\ldots,N-2\,,
\label{intforimm}
\end{eqnarray}
where $P_k$ are rank-1 orthogonal projectors which satisfy the 
completely integrable $2$-D 
$\field{C}P^{N-1}$ sigma model 
\begin{eqnarray}
\partial [\bar{\partial} P_k, P_k] + \bar{\partial} [\partial P_k, P_k]=0\,,
\qquad
P_k^2=P_k\,, \qquad 
P_k^{\dagger}=P_k\,.
\label{cpnmodel}
\end{eqnarray}

We first demonstrate that for any solution of the 
$\field{C}P^{N-1}$ sigma model (\ref{cpnmodel}) defined on the 
Riemann sphere $S^2$ with a finite action functional, the generalized 
Weierstrass formula for immersion of $2$-D surfaces (\ref{intforimm}) 
can be integrated explicitly up to a constant of integration and expressed 
in terms of the projectors $P_k$
\begin{eqnarray}
X_k(\xi, \bar{\xi}) = -i(P_k + 2 \sum_{j=0}^{k-1} P_j)\,,
\qquad
k=0,1,\ldots,N-2\,.
\label{eqforsurweierstrass} 
\end{eqnarray}

Indeed if we assume that the $\field{C}P^{N-1}$ sigma model is defined 
on the sphere $S^2$ with a finite action functional, 
then the complete set of regular solutions are known 
\cite{DinZakrzewski, Sasaki}. As a result, one gets three classes of 
solutions, namely (i) holomorphic (i.e. $\bar{\partial} f = 0$), 
(ii) antiholomorphic (i.e. ${\partial} f = 0$) 
and (iii) mixed. The mixed solutions can be determined from either 
the holomorphic or the antiholomorphic nonconstant functions by the 
successive application of the operator $P_{+}$ \cite{Zakrzewski} 
\begin{eqnarray}
P_{+}: f \in \field{C}^N \rightarrow
P_{+}f=\partial f -f \frac{f^{\dagger} \partial f}{f^{\dagger} f}\,,
\label{operatorintermsoff}
\end{eqnarray}
where $f$ is any nonconstant holomorphic function. This allows one to 
construct mixed solutions $f_k = P_+^kf$ which represent harmonic maps 
from $S^2$ to the $\field{C}P^{N-1}$ model. Here, the operator $P_+^k$ 
is obtained by applying successively $k$ times the operator $P_+$ 
and $P_+^0 =$ id. Hence, using (\ref{operatorintermsoff}) for 
every $k \leq N-1$ we can define a set of rank-1 projectors 
$\{P_0, P_1, \ldots, P_k\}$ 
\begin{eqnarray}
P_k: = \frac{f_k \otimes f_k^{\dagger}}{f_k^{\dagger} \cdot f_k}\,,
\qquad
k = 0, 1, \ldots. N-1\,,
\label{seqprojectors}
\end{eqnarray}
which determine conservation laws 
of the form (\ref{cpnmodel}). The first ($k=0$) and the last 
($k=N-1$) conservation laws are related to the holomorphic and 
antiholomorphic solutions, respectively, while 
the intermediate ones are related to the mixed solutions. 
Consequently, according to the Weierstrass procedure 
we can obtain $2$-D surfaces for each projector 
$P_k$. By straightforward calculation one gets \cite{gy} 
\begin{eqnarray}
&&[\partial P_k, P_k] = 
\partial P_k + 2 
\frac{({P}_+^{k} f) \otimes ({P}_+^{k-1} f)^{\dagger}}{|{P}_+^{k-1} f|^2}\,, 
\nonumber \\
&&[\bar{\partial} P_k, P_k] = 
-\bar{\partial} P_k - 2 
\frac{({P}_+^{k-1} f) \otimes ({P}_+^{k} f)^{\dagger}}{|{P}_+^{k-1} f|^2}\,.
\label{commutatorwithprank1bar}
\end{eqnarray}
Hence, for every $k\leq N-1$, the Weierstrass formula for immersion 
(\ref{intforimm}) takes the form
\begin{eqnarray}
dX_k \!=\! -i\! \left[
(\partial P_k + 2 
\frac{({P}_+^{k} f) \otimes ({P}_+^{k-1} f)^{\dagger}}{|{P}_+^{k-1} f|^2}) d\xi + 
(\bar{\partial} P_k + 2 
\frac{({P}_+^{k-1} f) \otimes ({P}_+^{k} f)^{\dagger}}{|{P}_+^{k-1} f|^2}) d\bar{\xi}
\right]\!\!.
\label{weierstrassdataforall}
\end{eqnarray}
Note that the two surfaces corresponding to $k=1$ and $k=N-1$ are 
precisely the same objects, since one gets antiholomorphic solutions 
of the model after applying $N-1$ times the operator $P_+$ to the 
nonconstant holomorphic function $f$. Hence, there appear at most 
$N-2$ different surfaces as a result of the 
generalized Weierstrass formula. 

The integration of (\ref{weierstrassdataforall}) gives us 
(\ref{eqforsurweierstrass}) which 
can be shown as follows: It is immediately seen that 
for $k=0$ we have $X_0=-iP_0$ and upon differentiation we obtain 
\begin{eqnarray}
dX_0 \!=\! -i\! \left[
\partial P_0  d\xi + 
\bar{\partial} P_0  d\bar{\xi}
\right],
\label{weierstrassdataforall2}
\end{eqnarray}
which coincides with (\ref{weierstrassdataforall}) for $k=0$. For 
$k=1$, we need to show that 
\begin{eqnarray}
\partial P_0 = \frac{(P_+f)\otimes f^{\dagger}}{|f|^2}\,,
\end{eqnarray}
which could easily be computed by differentiating $P_0$ and 
bearing in mind that $f$ is holomorphic. In order to show that 
(\ref{eqforsurweierstrass}) holds for any $k$ we assume that 
\begin{eqnarray}
\partial(P_0 + P_1 + \cdots, P_{k-2}) = \frac{(P_+^{k-1}f)\otimes (P_+^{k-2}f)^{\dagger}}{|P_+^{k-2}f|^2}\,,
\label{assumption}
\end{eqnarray}
and then compute $\partial P_{k-1}$
\begin{eqnarray}
\partial P_{k-1} &=& \partial \left[
\frac{(P_+^{k-1}f)\otimes (P_+^{k-1}f)^{\dagger}}{|P_+^{k-1}f|^2}
\right]
\nonumber \\
&=& \frac{(P_+^{k}f)\otimes (P_+^{k-1}f)^{\dagger}}{|P_+^{k-1}f|^2} - 
\frac{(P_+^{k-1}f)\otimes (P_+^{k-2}f)^{\dagger}}{|P_+^{k-2}f|^2}\,,
\label{inductionprocess}
\end{eqnarray}
where we have used the fact that 
\begin{eqnarray}
\partial (P_+^{k-1}f)^{\dagger} = -\frac{|P_+^{k-1}f|^2}{|P_+^{k-2}f|^2}(P_+^{k-2}f)^{\dagger}\,,
\end{eqnarray}
together with the orthogonality relation
\begin{eqnarray}
(P_+^{k}f)^{\dagger}\cdot (P_+^{l}f) = 0\,, \qquad {\rm for} \quad k\neq l\,.
\end{eqnarray}
Thus, we have shown that 
\begin{eqnarray}
\partial \left(
\sum_{j=0}^{k-1} P_j
\right) = \frac{(P_+^{k}f)\otimes (P_+^{k-1}f)^{\dagger}}{|P_+^{k-1}f|^2}\,,
\end{eqnarray}
which indeed justifies (\ref{eqforsurweierstrass}). 

It is interesting to note that this immersion function $X_k$ 
coincides with the results obtained in \cite{gy}, namely 
the surface for nonholomorphic Veronese type solutions lives in 
$\field{R}^3$. 

For the $\field{C}P^{N-1}$ sigma model, it may also be of interest to 
investigate the links between the generalized Weierstrass formula for the 
immersion of $2$-D surfaces in the $su(N)$ algebra and the approach 
based on the linear spectral problem for constructing infinitely many 
surfaces in multi-dimensional Euclidean spaces $\field{R}^{N^2-1} \simeq su(N)$. 

Following the approach proposed by A. Sym and J. Tafel in 
\cite{Sym1, Sym2, Tafel, Sym3}, in particular using their formula 
\begin{eqnarray}
X_k(\xi, \bar{\xi}) = \alpha(\lambda)\phi_k^{-1}\frac{\partial \phi_k}{\partial \lambda}\,,
\qquad
\phi_k \in SU(N)\,,
\quad
k=0,1,\ldots,N-2\,,
\label{symtafel}
\end{eqnarray}
for integrable surfaces derived via the Lax pair \cite{Mikhailov, Zakharov} 
\begin{eqnarray}
\partial \phi_k = \frac{2}{1+\lambda}[\partial P_k, P_k] \phi_k,
\quad
\bar{\partial} \phi_k = \frac{2}{1-\lambda}[\bar{\partial} P_k, P_k] \phi_k, 
\quad
k=0, 1, \ldots,N-1,
\label{linearscattprob}
\end{eqnarray}
we demonstrate that there exist $2$-D surfaces with $su(N)$-valued 
immersion functions $X_k(\xi,\bar{\xi})$ which are precisely of the form 
(\ref{eqforsurweierstrass}). Here, $\lambda$ is a spectral parameter and 
the compatibility conditions for the system (\ref{linearscattprob}) 
coincide with the $\field{C}P^{N-1}$ model equations (\ref{cpnmodel}). 
Furthermore, $\alpha(\lambda)$ is some scalar function of $\lambda$. 

In the purely instantonic case (i.e. holomorphic and 
antiholomorphic solutions), the orthogonal projector $P$ has the form 
\begin{eqnarray}
P_0 = \frac{f \otimes f^{\dagger}}{f^{\dagger}\cdot f}\,,
\label{proforp0forproof}
\end{eqnarray}
which satisfies 
\begin{eqnarray}
[\partial P_0, P_0] = \frac{P_+f\otimes f^{\dagger}}{f^{\dagger}\cdot f}\,,
\qquad
[\bar{\partial} P_0, P_0] = -\frac{f\otimes (P_+f)^{\dagger}}{f^{\dagger}\cdot f}\,.
\end{eqnarray}
Looking for a solution $\phi_0 = \phi_0(\lambda)$ of the linear spectral problem 
(\ref{linearscattprob}) when $\phi_0$ tends to $1$ as 
$\lambda \rightarrow \infty$, we make the Ansatz \cite{Zakrzewski} 
\begin{eqnarray}
\phi_0 = I_{N} - \frac{2}{1-\lambda}P_0\,,
\label{phiforproof}
\end{eqnarray}
where $I_N$ is the $N \times N$ identity matrix. 
The inverse matrix of $\phi_0$ is given by 
\begin{eqnarray}
\phi_0^{-1} = I_{N} - \frac{2}{1+\lambda}P_0\,.
\label{phiinvforproof}
\end{eqnarray}
Hence, according to the Sym-Tafel formula (\ref{symtafel}), the 
surface associated with the $\field{C}P^{N-1}$ model is given 
up to an additive $su(N)$-valued constant by 
\begin{eqnarray}
X_0 = \frac{2}{1-\lambda^2} P_0\,.
\label{surforins}
\end{eqnarray}

For the nonholomorphic solutions (i.e. the mixed solutions), we proceed 
in an analogous way and find that \cite{Zakrzewski}\,,
\begin{eqnarray}
\phi_k = I_N + \frac{4\lambda}{(1-\lambda)^2} \sum_{j=0}^{k-1}P_j - 
\frac{2}{1-\lambda}P_k\,,
\qquad
k=1,\ldots,N-2\,,
\label{phigenforproof}
\end{eqnarray}
and the inverse of $\phi_k$ has the form
\begin{eqnarray}
\phi_k^{-1} = I_N - \frac{4\lambda}{(1+\lambda)^2} \sum_{j=0}^{k-1}P_j - 
\frac{2}{1+\lambda}P_k\,.
\label{phiinvgenforproof}
\end{eqnarray}
The Sym-Tafel formula (\ref{symtafel}) for the immersion function 
$X_k(\xi,\bar{\xi})$ of $2$-D surfaces associated with the 
$\field{C}P^{N-1}$ model is given up to an additive $su(N)$-valued 
constant by the formula
\begin{eqnarray}
X_k(\xi, \bar{\xi}) = \frac{2}{1-\lambda^2}\left(P_k + 2\sum_{j=0}^{k-1}P_j 
\right)\,,
\qquad
k=0,1,\ldots,N-2\,.
\label{finalsymtafelres}
\end{eqnarray}
This result for immersion of $2$-D surfaces in the $su(N)$ Lie 
algebras coincides with the one obtained from the Weierstrass 
representation (\ref{eqforsurweierstrass}) and shows the equivalence 
between these two approaches. 

To conclude, we give an explicit expression for $2$-D surfaces, 
obtained from the generalized Weierstrass formula and demonstrate 
that these surfaces coincide with the ones obtained from 
the Sym-Tafel formula.

\section*{Acknowledgments}
This work is supported in part by research grants from NSERC of Canada. 
\.{I}.Y. acknowledges a postdoctoral fellowship awarded by the Laboratory 
of Mathematical Physics of the CRM, Universit\'{e} de Montr\'{e}al.

\end{document}